\documentclass[psamsfonts]{amsart}
\usepackage[utf8]{inputenc}
\usepackage[english]{babel}
\usepackage[T1]{fontenc}
\usepackage{amsfonts}
\usepackage{amsmath}
\usepackage{amsfonts}
\usepackage{bbm}
\usepackage{graphicx}
\usepackage{array}
\usepackage{lipsum}
\usepackage{amsaddr}

\usepackage{adjustbox}
\usepackage{transparent}

\usepackage{isomath}
\usepackage{mathtools}
\usepackage{amsthm}
\usepackage{slashed}
\usepackage{amssymb}
\usepackage{enumitem}
\usepackage{hyperref}
\hypersetup{colorlinks=true, 
citecolor=red}
\usepackage[new]{old-arrows}

\newtheorem{thm}{Theorem}[section]
\newtheorem{cor}[thm]{Corollary}
\newtheorem{prop}[thm]{Proposition}
\newtheorem{lem}[thm]{Lemma}

\theoremstyle{definition}
\newtheorem{defn}[thm]{Definition}

\newtheorem{rmk}[thm]{Remark}

\newtheorem{exmp}[thm]{Example}

\theoremstyle{remark}

\usepackage{tikz-cd} 
\usetikzlibrary{positioning, shapes, arrows.meta}
\tikzset{
    answer/.style={rectangle, draw, text width=15em, text badly centered, node distance=1cm, inner sep=0pt, minimum height=4em},
    block/.style={rectangle, draw, text width=10em, text centered},
     block2/.style={rectangle, draw, text width=5em, text centered},
     block3/.style={rectangle, draw, text width=5em, text centered, color=white},
}

\DeclareMathSymbol{\upLambda}{\mathalpha}{operators}{3}

\title[ A Berezin-Toeplitz Quantization viewpoint]{Structure preserving discretization:\\ A Berezin-Toeplitz Quantization viewpoint}
\author{Damien Tageddine, Jean-Christophe Nave}
\address{Department of Mathematics and Statistics, McGill University}
\begin{document}
\maketitle
\begin{abstract}
In this paper, we introduce a comprehensive axiomatization of structure-preserving discretization through the framework of commutative diagrams. By establishing a formal language that captures the essential properties of discretization processes, we provide a rigorous foundation for analyzing how various structures—such as algebraic, geometric, and topological features—are maintained during the transition from continuous to discrete settings. Specifically, we establish that the transition from continuous to discrete differential settings invariably leads to noncommutative structures, reinforcing previous observations on the interplay between discretization and noncommutativity.\\
We demonstrate the applicability of our axiomatization by applying it to the Berezin-Toeplitz quantization, showing that this quantization method adheres to our proposed criteria for structure-preserving discretization. We establish in this setting a precise limit theorem for the approximation of the Laplacian by a sequence of matrix approximations. This work enriches the theoretical understanding of discretization and sets the stage for its broader applications to various discretization methods.
\end{abstract}
\newpage
\section{Introduction}
\noindent
The discretization of continuous structures plays a pivotal role in both theoretical and applied mathematics, particularly in fields such as numerical analysis, computational physics, and data science. Traditional approaches to discretization often rely on heuristic or ad hoc methods, which can lead to inconsistencies and inefficiencies. Recent advances suggest that a more rigorous approach to discretization lies in preserving the essential features of the continuous structures while transitioning to a discrete framework. A key aspect of effective discretization is ensuring that the discrete model retains the underlying geometric and algebraic properties of the continuous one.\\

\noindent
Structure-preserving discretization methods have gained significant attention due to their ability to retain key properties of continuous systems, such as symmetries, conservation laws, and geometric structures, when transitioning to discrete models. State-of-the-art approaches include finite element exterior calculus (FEEC) \cite{christiansen_construction_2008, arnold_finite_2010}, which ensures stability and convergence in numerical approximations of partial differential equations by preserving the de Rham complex. Geometric integration methods, such as symplectic integrators, maintain the symplectic structure in Hamiltonian systems, ensuring long-term accuracy in energy conservation \cite{christiansen_topics_2011,gawlik_geometric_2011}. Discrete exterior calculus (DEC) extends these ideas to computational geometry by discretizing differential forms while preserving topological properties like divergence and curl \cite{hirani_discrete_2003, desbrun_discrete_2005}. Recent advances have also focused on the relation between noncommutative geometries and discrete differential structures \cite{lundervold_hopf_2010, mclachlan_butcher_2017}. This is in particular exemplified by the Berezin-Toeplitz quantization. This quantization discretizes certain types of Poisson algebras while maintaining algebraic structures in the transition from continuous to discrete systems \cite{modin_eulerian_2023}. Hence, quantization provides a robust framework for discrete geometry applications.\\

\noindent
Although there are similarities among the different approaches of structure-preserving discretization, it is yet not clear how one can describe under a unified theory, eclectic methods ranging from finite exterior calculus to quantization; we identify this as a gap in the theory.\\

\noindent
In this paper, we present an axiomatization of structure-preserving discretization by leveraging the language of commutative diagrams. Our primary objective is to formalize the process of discretization such that it maintains the integrity of underlying geometric and algebraic structures. We argue that by employing commutative diagrams, one can systematically capture and enforce the preservation of essential relationships between continuous and discrete models. A subsequent objective is to demonstrate that the structure-preserving discretization of differential structures, following the previous axioms, lead to noncommutative differential structure on operator algebras.\\

\noindent
Commutative diagrams provide a visual and conceptual tool from category theory that allow us to represent the interrelations between various mathematical structures and their discretizations. By constructing and analyzing these diagrams, we can ensure that the discretization process respects the morphisms and relationships inherent in the continuous setting. This approach not only facilitates a clearer understanding of the discretization process but also helps in developing algorithms that are both theoretically sound and practically effective.\\

\noindent
We demonstrate, by applying our axiomatization, that the discretization of a continuous differential structure is inherently given by a noncommutative 
geometry. This means that the resulting discrete structure is represented by operator algebras and the differentiation is realized by a commutator with a selfadjoint operator $D$. The results presented in this work corroborate previous observations and constructions \cite{tageddine_noncommutative_2022, tageddine_statistical_2023}.
The dictionary between continuous and discrete settings is identical to the dictionary between classical geometry and noncommutative one. This dictionary is summarized in the following table.
\begin{table}[h!]
\center
\parbox{1\linewidth}{
\centering
\begin{tabular}{|c|c|c|}
\hline 
Structure & Continuous & Discrete \\ 
\hline 
Space & Smooth manifold $M$ & simplex or graph \\ 
\hline 
Algebra & Algebra of functions $C^\infty(M)$ & Operator algebra $ \mathfrak{B}$ \\ 
\hline 
Differential & Differential form $df $ & Commutator  $[D,F]$ \\ 
\hline 
Integration & $\int f(x)dx$ & $\mathrm{Tr}_{\omega}(F)$ \\
\hline 
\end{tabular} 
\caption{Noncommutative geometry: from continuous to discrete}
\label{tab1}
}
\end{table}

\noindent
The Berezin-Toeplitz quantization is an archetypal example of deformation of a classical geometry given by a compact Kähler manifold $(M,g,\lbrace\cdot,\cdot\rbrace)$ to a noncommutative finite dimensional operator algebra $(M_n(\mathbb{C}),\left[ \cdot, \cdot \right] )$. Thus, it serves as a fundamental example of transition from a continuous to a discrete space. We use it in this work as a central application of our axiomatization in order to show that it satisfies the conditions of a structure preserving discretization. \\

\noindent
The paper is organized as follows: Section \ref{Sect1} introduces the key concepts and axioms that form the foundation of our framework. Section \ref{Sect2} demonstrates how commutative diagrams can effectively formalize the preservation of structures and explains how noncommutative geometry emerges from the structure-preserving discretization of a differential algebra. Section \ref{Sect3} provides a comprehensive review of diverse examples of structure-preserving discretizations, showcasing the framework’s flexibility. Finally, Sections \ref{Sect4} and \ref{Sect5} apply these findings to the Berezin-Toeplitz quantization, highlighting its adherence to the proposed criteria; Section \ref{Sect5} being a walk through example on the sphere.\\
Through this work, we hope to contribute to a deeper understanding of discretization methods and to provide a foundation for future research in this area.
\section{Axiomatization of the theory of discretizations}
\label{Sect1}
\noindent
In the first section, we introduce our axiomatization of structure-preserving discretization. The axioms rely on Category theory tools that we introduce briefly.  
\subsection{Category theory}
Let us recall that a category $\mathcal{C}$ is defined by a class $\mathrm{ob}(\mathcal{C})$ whose elements are called \textit{objects} and a class $\mathrm{hom}(\mathcal{C})$ whose elements are called \textit{morphisms} or \textit{arrows}. Each arrow $f$ has a \textit{source} object $X$ and a \textit{target} object $Y$. A morphism $f$ from $X$ to $Y$ is then denoted $f:X\rightarrow Y$. The class of arrows $\mathrm{hom}(\mathcal{C})$ is equipped with a binary operation called \textit{composition of arrows} between two morphism $f$ and $g$:
\begin{equation*}
f:X\rightarrow Y, \quad g:Y\rightarrow Z \quad \text{define}\quad g\circ f:X\rightarrow Z
\end{equation*}
such that $f\circ g$ is an element of $\mathrm{hom}(\mathcal{C})$. The composition is associative and there is exactly one \textit{identity morphism} denoted $1_X$ for every object $X$ such that for every arrow $f:X\rightarrow Y$:
\begin{equation*}
1_Y\circ f =f =f\circ f_X.
\end{equation*}
One says that a morphism $f:X\rightarrow Y$ is invertible and called an \textit{isomorphism} (or more simply an invertible arrow) if there exists a morphism $g:Y\rightarrow X$ such that $f\circ g=1_Y$ and $g\circ f=1_X$.
\subsection{Axiomatization} In what will follow, we will always, unless stated otherwise, consider objects in the category of Banach spaces. Nevertheless, one should have in mind that many of the constructions can be extended to merely topological spaces.
\begin{defn} 
\label{def1}
A \textit{discretization} $\mathfrak{D}(f)$ of an arrow $(f:C_1\rightarrow C_2)$ is a sequence of arrows $(f_n:C_1^n\rightarrow C_2^n)_{n\in \mathbb{N}}$ producing the following diagram
\begin{equation}
\begin{tikzcd}
 C_1 \arrow[r,"f"] \arrow[d,"\pi^n_1"'] & C_2  \arrow[d,"\pi^n_2"] \\
C^n_1 \arrow[r,"f_n"] & C^n_2 
\end{tikzcd}
\label{diagram}
\end{equation}
for each $n\in \mathbb{N}$. The maps $\pi^n_i:C_i\rightarrow C^n_i$, $i=1,2$ are surjective and contractive linear maps for all $n\in \mathbb{N}$. In addition, these maps satisfy the limit conditions:
\begin{equation}
\lim_{n\rightarrow \infty}\|\pi_i^nx\|_{C^n_i}=\|x\|_{C_i}, \quad \text{for all $x\in C_i$ where $i=1,2$.}
\end{equation}
A discretization $\mathfrak{D}(f)$ is \textit{faithful}   if:
\begin{equation*}
(\text{$f$ is an invertible arrow) $\Longrightarrow$ ($f_n\in \mathfrak{D}(f)$ is invertible for all $n\in \mathbb{N}$}).
\end{equation*}
\noindent
In the special case where the discretization $\mathfrak{D}(f:C\rightarrow C)$ does not depend on the map $f$, i.e. the maps $\pi^n$ are independent of $f$, then we say that we have a discretization of the object $C$ and we write $\mathfrak{D}(C)$.
\end{defn}
\begin{rmk}
The previous definitions generalize mutatis mutandis to the discretization $\mathfrak{D}(\lbrace f_n\rbrace_{n\in \mathbb{N}})$ of a sequence of arrows $\lbrace f_n\rbrace_{n\in \mathbb{N}}$ in the category $\mathcal{C}$. 
\end{rmk}
\begin{defn}[Structure preservation]
\label{def2}
Let $C_1$ and $C_2$ be objects in a category $\mathcal{C}$. Consider $f\in \mathrm{hom}(C_1,C_2)$ an arrow in $\mathcal{C}$ such that $(C_1,f)$ defines an object in a category $\mathcal{B}$. A discretization $\mathfrak{D}(f)$ is \textit{structure preserving} if it satisfies the following conditions:
\begin{itemize}
\item[1)] $(C^n_1,f_n)\in \mathrm{ob}(\mathcal{B})$ and $C_1^n\in \mathrm{ob}(\mathcal{C})$ for all $n\in \mathbb{N}$.
\item[2)] the diagram \eqref{diagram} commutes asymptotically:
\begin{equation*}
\|f_n\circ\pi_1^n(x)-\pi_2^n\circ f(x)\|\longrightarrow 0 \quad \text{as $n\rightarrow \infty$ and for all $x\in C_1$}
\end{equation*}
In such case, we will say that the discretization $\mathfrak{D}(f)$ is \textit{consistent}. 
\end{itemize}
\end{defn}
\begin{exmp}[Euler method]
Consider the algebra $C^\infty(\mathbb{S}^1)$ of smooth functions over the circle. For every $n\in \mathbb{N}$, let $X_n=\lbrace x_1,x_2,\dots,x_n\rbrace$ be a finite collection of points on $\mathbb{S}^1$ such that
$d(x_i,x_{i+1})=\frac{2\pi}{n},$ for $ i=1,2,\dots,n-1$ and where $d$ is the Euclidean distance on the circle.
Consider now the right-shift map on $\mathbb{C}^n$
\begin{equation*}
S:\mathbb{C}^n\rightarrow \mathbb{C}^n \quad S(y_1,y_2,\dots,y_n)=(y_n,y_1,y_2,\dots,y_{n-1})
\end{equation*}
and define the family Euler operators by 
\begin{equation*}
E_n=\frac{n}{2\pi}(S-1)
\end{equation*}
Moreover, we define the maps
\begin{equation*}
\pi^n:C^\infty(\mathbb{S}^1)\rightarrow \mathrm{Hom}(X_n,\mathbb{C}^n)\qquad \pi^n(f)=(f(x_1),f(x_2),\cdots,f(x_n)).
\end{equation*}
Consider now the differential operator 
\begin{equation*}
\frac{d}{d\theta}:C^\infty(\mathbb{S}^1)\rightarrow C^\infty(\mathbb{S}^1) \qquad f\mapsto \frac{df}{d\theta}
\end{equation*}
such that $(C^\infty(\mathbb{S}^1),\frac{d}{d\theta})$ is an object of the category of differential algebra.\\

\noindent
We can now look at the Euler discretization method of the pair $(C^\infty(\mathbb{S}^1),\frac{d}{d\theta})$ defined by the following diagram:
\begin{equation*}
\begin{tikzcd}
C^\infty(\mathbb{S}^1) \arrow[r,"\frac{d}{d\theta}"] \arrow[d,"\pi^n"'] & C^\infty(\mathbb{S}^1) \arrow[d,"\pi^n"] \\
\mathrm{Hom}(X_n,\mathbb{C}^n) \arrow[r,"E_n"] & \mathrm{Hom}(X_n,\mathbb{C}^n) 
\end{tikzcd}
\end{equation*}
We readily verify 
\begin{equation*}
\lim_{n\rightarrow \infty}\|\pi^n(f)\|_\infty = \|f\|_\infty \quad \text{and} \quad \|E_n\circ\pi^n(f)-\pi^n\circ \frac{d}{d\theta}(f)\| \longrightarrow 0
\end{equation*}
and thus, the Euler method defines a discretization. Therefore, it would be tempting to conclude that the Euler discretization is a structure preserving discretization of the object $(C^\infty(\mathbb{S}^1),\frac{d}{d\theta})$. However, this is not the case since the pair $(\mathrm{Hom}(X_n,\mathbb{C}^n),E_n)$ does not define a differential algebra, since $E_n$ does not satisfy the Leibniz rule.
\end{exmp}
\begin{rmk}
In fact, a structure-preserving discretization of the differential algebra $(C^\infty(\mathbb{S}^1),\frac{d}{d\theta})$ necessitates the use of a Lie bracket  in a finite dimensional context, which indeed satisfies the product rule and defines a derivation. This will be explored in the next section.
\end{rmk}
\noindent
We can also require stronger conservation properties between the continuous object $C_1$ and its discretization $C_1^n$.
\begin{defn}[Strongly structure preserving]
Let $C_1$ and $C_2$ be objects in a category $\mathcal{C}$. Consider $f\in \mathrm{hom}(C_1,C_2)$ an arrow in $\mathcal{C}$ such that $(C_1,f)$ defines an object in a category $\mathcal{B}$. A discretization $\mathfrak{D}(f)$ is \textit{strongly structure preserving} if it satisfies the following conditions:
\begin{itemize}
\item[1)] $(C^n_1,f_n)\in \mathrm{ob}(\mathcal{B})$ and $C_1^n\in \mathrm{ob}(\mathcal{C})$ for all $n\in \mathbb{N}$.
\item[2)] the diagram \eqref{diagram} commutes for all $n\in \mathbb{N}$;
\item[3)] the maps $\pi_i^n:C_i\rightarrow C_i^n$ are surjective homomorphisms for all $n\in \mathbb{N}$ and for $i=1,2$ ( i.e. $\pi_i^n\in \mathrm{hom}(\mathcal{C})$ for all $n\in \mathbb{N}$).
\end{itemize}
\end{defn}
\begin{exmp} Let $\mathcal{C}$ the category of prehilbert space and $\mathcal{B}$ be the category of differential algebra. Then, for instance, the pair $(C^\infty(\mathbb{S}^1),\left\langle \cdot ,\cdot \right\rangle )$ consisting of the space smooth function on the circle with its usual scalar product, is an object of $\mathcal{C}$; and the pair $(C^\infty(\mathbb{S}^1),\frac{d}{d\theta})$ defining a differential algebra is an object of $\mathcal{B}$.\\
A strongly structure preserving discretization of the differential algebra $(C^\infty(\mathbb{S}^1),\frac{d}{d\theta})$ (if it exists) is a commutative diagram:
\begin{equation*}
\begin{tikzcd}
(C^\infty(\mathbb{S}^1),\left\langle \cdot ,\cdot \right\rangle ) \arrow[r,"\frac{d}{d\theta}"] \arrow[d,"\pi^n"'] & (C^\infty(\mathbb{S}^1),\left\langle \cdot ,\cdot \right\rangle ) \arrow[d,"\pi^n"] \\
(\mathcal{A}_n ,\left\langle \cdot ,\cdot \right\rangle_n )\arrow[r,"d_n"] & (\mathcal{A}_n ,\left\langle \cdot ,\cdot \right\rangle_n )
\end{tikzcd}
\end{equation*}
such that:
\begin{itemize}
\item[1)] $(\mathcal{A}_n,d_n)$ is a differential algebra;
\item[2)] $\frac{d}{d\theta}\circ \pi^n = \pi^n\circ d_n$ for all $n$;
\item[3)] $(\mathcal{A}_n ,\left\langle \cdot ,\cdot \right\rangle_n )$ is a prehilbert space and $\pi^n$ a partial isometry (i.e. an isometry on its support) that is:
\begin{equation*}
\left\langle \pi^nx,\pi^ny\right\rangle_n = \left\langle x,y\right\rangle
\end{equation*}
for $x,y$ in the support.
\end{itemize}
\end{exmp}
\subsection{Convergence} Consider a discretization $\mathfrak{D}(f)$ of an arrow $f$. Heuristically, the convergence of the discretization means that $(C_1^n,f_n)$ converges to $(C_1,f)$. Hence, to properly define a notion of convergence, one needs to compare elements between different levels of dicretization.\\

\noindent
Let us recall that, the graph of an arrow $f$, denoted $\mathrm{Gr}(f)$, is defined as the subset of $C_1\times C_2$ such that:
\begin{equation*}
\mathrm{Gr}(f)=\left\lbrace (x,y)\in C_1\times C_2 : y=f(x) \right\rbrace 
\end{equation*}
We denote by $p_i:\mathrm{Gr}(f)\rightarrow C_i$ for $i=1,2$ the obvious coordinate projections.
\begin{defn}[Convergence]
\label{conv}
Consider a discretization $\mathfrak{D}(f)$ of an arrow $f$. We say that $\mathfrak{D}(f)$ is \textit{convergent} if to any of the projection maps $\pi^n_i:C_i\rightarrow C^n_i$, we can associate an injective contractive linear map $s_i^n:C_i^n\rightarrow C_i$ such that
\begin{equation*}
\lim_{n\rightarrow +\infty}\|x-s_i^n\circ \pi^n_i(x)\|=0, \quad \text{for all $x\in p_i(\mathrm{Gr}(f))$, and for $i=1,2$.}
\end{equation*}
The maps $s_i^n$ will be called \textit{section} map.
\end{defn}
\noindent
For a structure preserving discretization $\mathfrak{D}(f)$, one can obtain an equivalent definition of convergence.
\begin{prop}
Consider a structure preserving discretization $\mathfrak{D}(f)$ of an arrow $f$. The discretization $\mathfrak{D}(f)$ is convergent if and only if to any of the projection maps $\pi^n_i:C_i\rightarrow C^n_i$, we can associate an injective contractive linear map $s_i^n:C_i^n\rightarrow C_i$ such that
\begin{equation*}
\lim_{n\rightarrow +\infty}\|x-s_1^n\circ \pi^n_1(x)\|=0, \qquad \text{and} \qquad
\lim_{n\rightarrow +\infty}\|f(x)-s_2^n\circ f_n (\pi^n_1(x))\|=0
\end{equation*}
\end{prop}
\begin{proof}
We prove the \textit{if} part of the proof, the \textit{only if} can be proven identically. Assume that section maps exist, then one get the following inequality
\begin{equation*}
\|f(x)-s_{2}^n\circ \pi_{2}^n(f(x)))\|\leq\|f(x)-s_{2}^n\circ f_n(\pi_1 ^n(x))\| +\|s_2^n(\pi_2^n\circ f(x)- f_n\circ\pi_1^n(x))\|
\end{equation*}
where the first term in the right-hand-side goes to zero by assumption and the second term also vanishes from the structure preserving condition.
\end{proof}
\noindent
If we assume that the arrow $f$ is also invertible, then one can deduce the following corollary relating convergence of $\mathfrak{D}(f)$ and $\mathfrak{D}(f^{-1})$.
\begin{cor}
Let $f$ be an invertible arrow, such that the discretizations $\mathfrak{D}(f)$ and $\mathfrak{D}(f^{-1})$ are structure preserving and faithful. Then, $\mathfrak{D}(f)$, respectively $\mathfrak{D}(f^{-1})$, converges if  
\begin{equation*}
\lim_{n\rightarrow +\infty}\|f(x)-s_2^n\circ f^{-1}_n (\pi^n_1(x))\|=0, \qquad \text{and} \qquad
\lim_{n\rightarrow +\infty}\|f(x)-s_2^n\circ f_n (\pi^n_1(x))\|=0
\end{equation*}
\end{cor}
\section{Derivations and Noncommutative geometry}
\label{Sect2}
\subsection{Derivations}
The central focus of discretization theory is the study of differential equations derived from the differential structure of a manifold. 
In the following paragraph, we aim to demonstrate how the structure-preserving discretization of a smooth differential operator inevitably results in a noncommutative differential structure.  \\

\noindent
 A differential structure means here the data of a unital algebra of smooth functions $A$ and a derivation; we recall that a derivation $d:A\rightarrow M$ is map such that:
\begin{itemize}
\item[1)] $d$ is linear,
\item[2)] $d(1)=0$,
\item[3)] $d(ab)=d(a)b + ad(b)$.
\end{itemize}
where $M$ is a $A$-bimodule. The data $(A,d)$ is called a \textit{differential algebra} and constitute one of the fundamental structure that one may want to preserve under discretization. Consider a structure preserving discretization $(A_n,d_n)$ of a given differential algebra $(A,d)$ given by the following diagram:
\begin{center}
\begin{tikzcd}
  A \arrow[r,"d"] \arrow[d,"\pi^n"']
    & M \arrow[d,"\pi^n"] \\
  A_n\arrow[r,"d_n"]
& M_n
\end{tikzcd}
\end{center}
\begin{prop} If $A_n$ is isomorphic to a matrix algebra then there exists a self-adjoint element $D_n$ such that
\begin{equation*}
d_n(a)=\left[ D_n,a\right] 
\end{equation*}
In addition, if the discretization $(A_n,d_n)$ is structure preserving, then
\begin{equation*}
 \lim_{n\rightarrow \infty}\|\pi^n\circ d(a)-\left[ D_n,\pi^n(a)\right] \| = 0
\end{equation*}
\end{prop}
\begin{proof}
The claim immediately follows from Theorem 4.1.6 in \cite[p. 156]{sakai_c-algebras_1998} and the axiom of a structure preserving discretization.
\end{proof}
\begin{rmk}
Structure preserving discretizations $(A_n,d_n)$ of differential algebra are also called \textit{approximately inner derivations} in the theory of Operator Algebras, see for instance \cite{sakai_operator_2008}. The last proposition provides a first bridge between application of $C^*$-algebra theory (or operator theory) and structure preserving discretizations.
\end{rmk}
We have established that one approximate an differential algebra $(A,d)$ of functions by a matrix algebra, then the discrete differential structure is inherently noncommutative. We now show a concrete case where such approximation exists.
\begin{defn}
A filtration of a Hilbert space $\mathcal{H}$ is a sequence $\mathcal{F}=\left\lbrace \mathcal{H}_1,\mathcal{H}_2,\dots \right\rbrace $ of finite dimensional subspaces of $\mathcal{H}$ such that $\mathcal{H}_n\subset\mathcal{H}_{n+1}$ and 
\begin{equation*}
\mathcal{H}=\overline{\bigcup_{n=1}^\infty \mathcal{H}_{n}}
\end{equation*}
Let $\mathcal{F}=\left\lbrace \mathcal{H}_n\right\rbrace $ be a filtration of $\mathcal{H}$ and let $P_n$ be the projection onto $\mathcal{H}_n$. The degree of an operator $D$ is defined by
\begin{equation*}
\mathrm{deg}(D)=\sup_{n\geq 1}\mathrm{rank}(P_nD-DP_n)
\end{equation*}
\end{defn}
\noindent
Therefore, one sees that if the domain $\mathcal{H}$ of a continuous differential operator $D$ admits a filtration family $\lbrace\mathcal{H}_n\rbrace$, then we can construct the following diagram of discretization:
\begin{equation}
\begin{tikzcd}
\mathcal{H} \arrow[r,"D"] \arrow[d,"P_n"'] & \mathcal{H}  \arrow[d,"P_n"] \\
\mathcal{H}_n \arrow[r,"D_n"] & \mathcal{H}_n 
\end{tikzcd}
\label{diagram2}
\end{equation}
Thus, the lack of commutativity of the previous diagram i.e. the existence of a structure preserving discretization, in this case is measured by the degree $\mathrm{deg}(D)$ of the operator $D$. We can then restrict ourselves to subclass of operators with vanishing degrees.
\begin{defn}[Block diagonal operator]
A linear operator $D$ on a Hilbert space $H$ is called \textit{block diagonal} if there exists an increasing sequence of finite rank projections, $P_1\leq P_2\leq P_3\leq \cdots$ such $\|[D,P_n]\|=\|DP_n-P_nD\|=0$ for all $n\in \mathbb{N}$ and $P_n\rightarrow 1_{\mathcal{H}}$ (in the strong operator topology) as $n\rightarrow \infty$.
\end{defn}
\noindent
Recall that two (possibly unbounded) self-adjoint operators $A$ and $B$ are said to \textit{commute} if and only if all the projections in their associated projection-valued measures commute. Hence, a trivial example of increasing sequence of commuting projections is given by spectral projectors. \\
Similarly to the axioms of structure preserving discretization, the block diagonal condition might be too restrictive; one can then relax the commutativity condition.
\begin{defn}
A linear operator $D$ on a Hilbert space $H$ is called \textit{quasidiagonal} if there exists an increasing sequence of finite rank projections, $P_1\leq P_2\leq P_3\leq \cdots$ such $\|[D,P_n]\|=\|DP_n-P_nD\|\rightarrow 0$ and $P_n\rightarrow 1_{\mathcal{H}}$ (in the strong operator topology) as $n\rightarrow \infty$.
\end{defn}
\noindent
We can now state a theorem on existence of structure preserving discretization.
\begin{thm}[First Existence Theorem]
Consider a differential algebra of functions $(A,d)$ such that $da=Da-aD$ such that $D$ is a quasidiagonal (not necessarily bounded) operator, then $(A,d)$ admits a structure preserving discretization. If in addition $D$ is block diagonal, then $(A,d)$ admits a strongly structure preserving discretization.
\end{thm}
\begin{proof}
Since $d$ is a quasidiagonal derivation in an algebra of functions $A$, there  exists a self-adjoint operator $D$ such that $da=[D,a]:=Da-aD$ and a family of finite rank projections such that $\|[D,p_n]\|\rightarrow 0$. Then, we can define the surjective contractions:
\begin{equation*}
\pi^n:A\rightarrow A_n, \qquad \pi^n(a)=p_nap_n
\end{equation*}
and thus satisfy
\begin{equation*}
\pi^n(da)=p_n[D,a]p_n=[p_nDp_n,\pi^n(a)]+(p_n[D,p_n]ap_n-p_na[D,p_n]p_n).
\end{equation*}
Therefore, we deduce the following upper-bound
\begin{equation*}
\|\pi^n(da)-[D_n,\pi^n(a)]\| \leq \|[D,p_n]\| \|a\|\rightarrow 0
\end{equation*}
which implies the structure-preserving condition.\\
If $D$ is block-diagonal, then it admits a family $\lbrace p_n\rbrace$ of finite rank projectors such that $p_nD=p_nD$. Therefore, one has 
\begin{equation}
\pi^n(da)=p_n[D,a]p_n=[D_n,\pi^n(a)]=d_n\pi_n(a).
\end{equation}
i.e. the discretization is strongly structure preserving.
\end{proof}
\begin{rmk}
From the Existence Theorem, we see that a structure preserving discretization of a differential algebra, when it exists, is given by the data of $(A_n,H_n,D_n)$ where $A_n$ is a matrix algebra, $H_n$ is its representation space and $D_n$ is a selfadjoint operator defining the differential structure. This is the general setting of noncommutative geometry \cite{connes_noncommutative_1994}.
\end{rmk}
\subsection{$C^*$-algebras as a general setting for structure preserving discretizations}
In fact, most cases of discretization and existence of projection maps $\pi^n$ and section maps $s^n$ fit it in the general framework of $C^*$-algebras. Let us briefly recall that a $C^*$-algebra is a Banach algebra equipped with an involution $*$, compatible with the norm: $\|x^*x\|=\|x\|^2$.
\begin{defn}[Nuclear $C^*$-algebra]
A $C^*$-algebra $A$ is said to be \textit{nuclear} if there exists a net of natural numbers $(n_\lambda)_{\mathrm{\Lambda}}$ and nets of contractive, completely positive maps $\varphi_\lambda :A\rightarrow M_{n_\lambda}(\mathbb{C})$ and $\psi_\lambda : M_{n_\lambda}(\mathbb{C}) \rightarrow A$ such that $\lim_{\mathrm{\Lambda}}\|a-\psi_\lambda(\varphi_\lambda(a))\|=0$ for every $a\in A$.
\end{defn}
\begin{thm}[Second Existence Theorem]
\label{2nd}
Consider a pair $(A,f)$ such that $A$ is also a nuclear $C^*$-algebra. There always exists a discretization $\mathfrak{D}(A,f)$ of $(A,f)$. Moreover, if $(A,f)$ is considered as an object of the category of linear spaces, then a strongly structure preserving discretization always exists.
\end{thm}
\begin{proof}
Since $A$ is a nuclear $C^*$-algebra, there exist discretization maps $\varphi_n:A\rightarrow M_{n}(\mathbb{C})$ and section maps $\psi_n:M_{n}(\mathbb{C}) \rightarrow A$. 
We need to verify that $\lim_{n\rightarrow \infty}\|\varphi(a)\|=\|a\|$; the limit follows from  the fact that $\lim_{n\rightarrow \infty}\|a-\psi_n\circ\phi_n(a)\|=0$ and the contraction property and:
\begin{equation}
\|\psi_n\circ\varphi_n(a)\|\leq \|\varphi_n(a)\|\leq \|a\|
\end{equation}
Therefore, by taking the limit we obtain that $\lim_{n\rightarrow\infty}\|\varphi_n(a)\|=\|a\|$. We can construct the following discretization diagram: 
\begin{equation*}
\begin{tikzcd}
A \arrow[r,"f"] \arrow[d,"\varphi^n"'] & A    \arrow[d,"\varphi^n"] \\
M_n(\mathbb{C}) \arrow[r,"f_n"] & M_n(\mathbb{C})
\end{tikzcd}
\end{equation*}
where the discrete arrow $(f_n)$ is defined by:
\begin{equation}
f_n(\pi^n(a))=\pi^n(f(a)).
\end{equation}
\end{proof}
\begin{rmk}
It turns out that all commutative $C^*$-algebra $C(X)$, of continuous functions on a compact space $X$, are nuclear. More generally, if $A$ and $B$ are nuclear $C^*$-algebras then $A\otimes B$ is a nuclear $C^*$-algebra. In particular, if $A=C(X)$ then $A\otimes B$ can be identified with $C(X,B)$ i.e. the $C^*$-algebra of continuous functions with values in $B$, under the map $(f\otimes b)(x)=f(x)b$.
\end{rmk}
\begin{exmp}
A fundamental example is given by the space of continuous differential forms on a manifold $X$ of dimension $n$. Recall that a continuous differential $r$-form is a map 
\begin{equation*}
\omega:X\rightarrow \mathrm{\Lambda}^r\mathbb{R}^n, \quad x\mapsto \sum_{|I|=r}\omega_I(x)dx^I
\end{equation*}
where the functions $\omega_I$ are continuous. If we identify the exterior algebra $\mathrm{\Lambda}\mathbb{R}^n$ with the Clifford algebra $\mathrm{Cl}(\mathbb{R}^n)$ then it acquires a $C^*$-algebra structre. The space of continuous differential forms is then isomorphic as a $C^*$-algebra to $C(X,\mathrm{Cl}(\mathbb{R}^n))$. Since the Clifford algebras are nuclear $C^*$-algebras, then the space of continuous differential forms is also a nuclear $C^*$-algebra and therefore admits a discretization according to the previous theorem.
\end{exmp}
\section{Zoology of structure preserving discretizations}
\label{Sect3}
\noindent
In the following section, we review archetypical examples of structure preserving discretizations through the lens of our axiomatization.
\subsection{Finite element exterior calculus}  We start by recalling some of the basic definitions of finite element exterior calculus.
\begin{defn}
Let $W^1,W^2,\dots$ be Hilbert spaces and let $d^k:W^k\rightarrow W^{k+1}$ be densely-defined closed linear operators satisfying
\begin{itemize}
\item[1)] $\mathrm{Im}(d^k)\subseteq \mathrm{Dom}(d^{k+1})$ and
\item[2)] $d^{k+1}d^k=0$
\end{itemize}
for all $k$. Then, we say that the data $\lbrace (W^k,d^k)\rbrace_{k}$ define a \textit{Hilbert complex}. The operators $d^k$ are called the \textit{differentials} of the complex. If all the differentials are bounded, the Hilbert complex is said to be bounded.
\end{defn}
\begin{defn}
A \textit{morphism of bounded Hilbert complexes} $\Phi:(V,d)\rightarrow (\tilde{V},\tilde{d})$ is a set of linear maps $\Phi^k:V^k\rightarrow \tilde{V}^k$ such that $\tilde{d}\Phi^k=\Phi^{k+1}d^k$ for all $k$. We indicate that $\Phi$ is a morphism by writing $\Phi \in \mathrm{hom}((V,d),(\tilde{V},\tilde{d}))$.
\end{defn}
\noindent
The framework of the finite element exterior calculus is given by the $L^2$-de Rham complex, represented by the following complex
\begin{equation*}
0\rightarrow \mathrm{H\Lambda}^0(\mathrm{\Omega})\xrightarrow[]{d}\mathrm{H\Lambda}^1(\mathrm{\Omega})\xrightarrow[]{d}\cdots \xrightarrow[]{d}\mathrm{H\Lambda}^n(\mathrm{\Omega})\rightarrow 0
\end{equation*}
where $\mathrm{\Omega}\subseteq \mathbb{R}^n$ is a Lipschitz domain. In order to give a numerical approximation of a PDE, the method builds a finite-dimensional subcomplex 
\begin{equation*}
0\rightarrow \mathrm{\Lambda}_\hbar^1(\mathrm{\Omega})\xrightarrow[]{d}\cdots \xrightarrow[]{d}\mathrm{\Lambda}_\hbar^n(\mathrm{\Omega})\rightarrow 0
\end{equation*} 
of the de Rham complex. In order to construct this subcomplex, one shows that there exists morphism $\pi^\hbar$ projecting the de Rham complex down to the appropriate subcomplex; so that each map $\pi^\hbar_k:\mathrm{H\Lambda}^k\rightarrow \mathrm{\Lambda}^k_\hbar$ is a projection onto the subspace $\mathrm{\Lambda}_\hbar^k$ and induces the following commuting diagram
\begin{equation}
\begin{tikzcd}
 H\mathrm{\Lambda}^k \arrow[r,"d"] \arrow[d,"\pi^\hbar_k"'] &   H\mathrm{\Lambda}^{k+1} \arrow[d,"\pi^\hbar_{k+1}"] \\
\mathrm{\Lambda}_\hbar^k \arrow[r,"d"] & \mathrm{\Lambda}_\hbar^{k+1} 
\end{tikzcd}
\end{equation}
In this case, the section maps $s_k^{\hbar}$ are simply identity maps. Therefore, if we consider $\lbrace H\mathrm{\Lambda}^k,d^k\rbrace$ as an object of the category of Hilbert complexes, then the data $(\lbrace \mathrm{\Lambda}^k_{\hbar},d\rbrace, \pi^\hbar_k)$ gives a strongly structure preserving discretization.
\begin{prop}
The Finite Element Exterior Calculus is a strongly structure preserving discretization.
\end{prop}
\begin{proof}
We verify that the discretization $(\lbrace \mathrm{\Lambda}^k_{\hbar},d\rbrace, \pi^\hbar_k)$ satisfies the axiom of a structure preserving discretization of the continuous data $\lbrace H\mathrm{\Lambda}^k,d^k\rbrace$ seen as an object of the category of Hilbert complexes. The following three conditions are fulfilled:
\begin{itemize}
\item[1)]$\lbrace \mathrm{\Lambda}^k_{\hbar},d\rbrace$ is a differential complex
\item[2)] the diagram commutes for all $n$:
\begin{equation*}
d\circ \pi^{\hbar}_k=\pi^{\hbar}_{k+1}\circ d
\end{equation*}
\item[3)] the maps $\pi^\hbar_{k}$ are bounded maps between normed spaces.
\end{itemize}
from which the statement follows.
\end{proof}
\subsection{Transfer Operator and diffeomorphism group}
In the following paragraph, we would like to illustrate how our axiomatization allows us to describe structure-preseving discretizations of diffeomorphism groups; we prove on practical examples how infinite group of symmetries are related to matrix groups.\\

\noindent
Consider a smooth compact manifold $X$ and its diffeomorphism group $\mathrm{Diff}(X)$. In this section, we identify structure preserving discretization of diffeomorphism groups which play a central role in PDE theory.\\

\noindent
The group $\mathrm{Diff}(X)$ admits a natural action on the Hilbert space $L^2(X)$, the space of square-summable half-densities on $X$; by choosing a smooth measure $m$ on $X$ one can identify $L^2(X)$ with the usual space of functions $f$ on $X$ which are square-summable with respect to $m$. Hence,  $\mathrm{Diff}(X)$ admits an obvious unitary representation on $H$ induced by its action on $X$: for any $\psi \in \mathrm{Diff}(X)$, define
\begin{equation*}
\widehat{\psi}:H\rightarrow H \qquad \widehat{\psi}(f)(x):= J_\psi(x)^{\frac{1}{2}}f(\psi^{-1}x)
\end{equation*}
where $J_\psi(x)=dm(\psi^{-1}x)/dm(x)$. We can restrict this representation to an operator $\widehat{\psi}:C^\infty(X)\rightarrow C^\infty(X)$ and study the possible structure preserving discretization of $\mathrm{GL}(C^\infty(X))$.\\

\noindent
Fix a diffeomorphism $\psi$, in order to present a structure preserving discretization of the pair $(C^\infty(X), \widehat{\psi})$, we need to specify in which category the object $(C^\infty(X), \widehat{\psi})$ belongs to: we will consider it as an object of the category of linear spaces. Therefore, following the definition, a structure preserving discretization is given by a commuting diagram
\begin{equation}
\begin{tikzcd}
C^\infty(X) \arrow[r,"\widehat{\psi}"] \arrow[d,"\pi^n"'] & C^\infty(X) \arrow[d,"\pi^n"] \\
A_N\arrow[r,"\widehat{\psi}_N"] & A_N
\end{tikzcd}
\label{diagram3}
\end{equation}
such that $\lim_{n\rightarrow \infty}\|\pi^n(f)\|_\infty = \|f\|_\infty $ and $(A_N,\widehat{\psi}_N)$ is a pair of a vector space $A_N$ and a linear map $\widehat{\psi}_N$.
\begin{lem}
There exists a family of surjective maps $\pi^n:C^\infty(X)\rightarrow \mathbb{C}^{N}$ such that $\lim_{n\rightarrow \infty}\|\pi^n(f)\|_\infty = \|f\|_\infty$ for all $f \in C^\infty(X)$.
\end{lem}
\begin{proof}
The algebra $C(X)$ of continuous functions is a nuclear abelian $C^*$-algebra, therefore it admits a discretization map $\pi^n:C(X)\rightarrow \mathbb{C}^{N}$; we restrict this map to $C^\infty(X)$. Moreover, the limit $\lim_{n\rightarrow \infty}\|\pi^n(f)\|=\|f\|$ follows from the Second Existence Theorem.
\end{proof}
\noindent
Therefore, diagram \eqref{diagram3} becomes 
\begin{equation}
\begin{tikzcd}
C^\infty(X) \arrow[r,"\widehat{\psi}"] \arrow[d,"\pi^n"'] & C^\infty(X) \arrow[d,"\pi^n"] \\
\mathbb{C}^N\arrow[r,"\widehat{\psi}_n"] & \mathbb{C}^N
\end{tikzcd}
\end{equation}
and the existence of $\widehat{\psi}_n$ is given by the Second Existence Theorem. We can then state the following theorem.
\begin{thm}
\label{thmdiff}
For any diffeomorphism $\psi\in \mathrm{Diff}(X)$, there exists a faithful, strongly structure preserving discretization of $\mathfrak{D}(C^\infty(X),\widehat{\psi})$.
\end{thm}
\begin{proof}
Following the proof of Theorem \eqref{2nd}, we define $\widehat{\psi}_n$ by
\begin{equation}
\widehat{\psi}_n\circ \pi^n(f)=\pi^n\circ \widehat{\psi}(f)
\label{identity}
\end{equation}
which gives us a linear surjection $\widehat{\psi}_n:\mathbb{C}^N\rightarrow \mathbb{C}^N$, since $\pi^n$ and $\widehat{\psi}$ are surjective; and thus $\widehat{\psi}_n$ is an isomorphism.
\end{proof}
\begin{rmk}
Since the previous discretization does not depend on the choice of diffeomorphism, we will abuse notation and denote by
\begin{equation*}
\mathfrak{D}(\mathrm{Diff}(X))=\mathrm{GL}_N(\mathbb{C}).
\end{equation*} 
to say that a discretization of the diffeomeorphism group $\mathrm{Diff}(X)$ is given by the finite dimensional group of invertible matrices $\mathrm{GL}_N(\mathbb{C})$.
\end{rmk}
\noindent
One could also be interested in subgroups of the diffeomorphism group. For instance, if the manifold $X$ is equipped with a volume form $d\omega$, we can define the scalar product between two function:
\begin{equation*}
\left\langle f,g \right\rangle =\int_X f\overline{g} \ d\omega
\end{equation*}
 Then, one can define the subgroup of diffeomorphisms preserving this inner product and denoted it by $\mathrm{Diff}^+_0(X)$.\\
 
 \noindent
One can also define the subgroup of markovian diffeomorphism, i.e. the group of diffeomorphism $\psi$ satisfying $\psi(1)=1$. We denote this subgroup by $\mathfrak{D}(\mathrm{Diff}_{\mathrm{m}}(X))$.\\

\noindent
Let us denote by $\mathrm{GL}^{\mathrm{st}}_N(\mathbb{C})$, the group of invertible stochastic matrices.
\begin{prop}
Let $(A_N,\widehat{\psi}_N)$ be a faithful strongly structure preserving discretization of the pair $(C^\infty(X), \psi)$ independent of the choice of $\psi$. If
\begin{itemize}
\item[i)] $C^\infty(X)$ is merely a vector space, then $\mathfrak{D}(\mathrm{Diff}(X))=\mathrm{GL}_N(\mathbb{C}).$
\item[ii)]$(C^\infty(X),\left\langle \cdot,\cdot \right\rangle )$ is a prehilbert space, then $\mathfrak{D}(\mathrm{Diff}_0^{+}(X))=\mathrm{SO}_N(\mathbb{C}).$
\item[iii)]$(C^\infty(X),1)$ is a vector space with a unit element (i.e. containing its scalar field), then $\mathfrak{D}(\mathrm{Diff}_{\mathrm{m}}(X))=\mathrm{GL}^{\mathrm{st}}_N(\mathbb{C}).$
\end{itemize}
\end{prop}
\begin{proof}
i) follows from Theorem \eqref{thmdiff}.\\
ii) If $(A_N,\widehat{\psi}_N)$ is a strongly structure preserving discretization of the pair $(C^\infty(X), \psi)$, then $(A_N, \left\langle \cdot , \cdot \right\rangle_N )$ has a prehilbert structure as well. Moreover, the map $\pi^N:C^\infty(X)\rightarrow A_N$ is an isometry on its support. Therefore, we have that for any $f,g\in C^\infty(X)$
\begin{equation}
\left\langle \widehat{\psi}_N\circ \pi^N(f), \widehat{\psi}_N\circ \pi^N(g) \right\rangle_N = \left\langle \pi^N\circ \widehat{\psi}(f) ,\pi^N\circ \widehat{\psi}(g) \right\rangle_N= \left\langle f,g \right\rangle .
\end{equation}
The statement follows from the surjectivity of $\widehat{\psi}$ and $\pi^n$.\\
iii) If $(A_N,\widehat{\psi}_N)$ is a strongly structure preserving discretization of the pair $(C^\infty(X), 1)$, then $(A_N, 1_N)$ has also a unit element. Moreover, the map $\pi^N:C^\infty(X)\rightarrow A_N$ is unital, i.e. $\pi^n(1)=1_N$. Using the commutativity axiom, we get, 
\begin{equation}
\widehat{\psi}_N(1_N)=\widehat{\psi}_N\circ \pi^N(1)=\pi^N\circ \widehat{\psi}(1)=1_N.
\end{equation}
Therefore, the discretization $\widehat{\psi}_N$ preserves the unit. Because  $\widehat{\psi}_N$ is invertible (from statement i)) then it is a left and right stochastic linear map.
\end{proof}
\noindent
In addition, we notice that a structure preserving diffeomorphism discretization of a diffeomorphism group induces a structure preserving discretization of the associated Lie algebra. Let us recall that the pullback representation gives the map
\begin{equation*}
\mathrm{Diff}(X)\rightarrow \mathrm{Lin}(C^\infty(X)) \quad \psi\mapsto (f\mapsto \widehat{\psi}(f))
\end{equation*} 
Differentiating it at the identity $\psi=\mathrm{id}$ gives a representation of the algebra of vector fields $\mathrm{Vect}(X)$ over the space $X$ into the space of derivations:
\begin{equation*}
\mathrm{Vect}(X)\rightarrow \mathrm{Der}(C^\infty(X)) \quad \Psi\mapsto \mathcal{L}_\Psi
\end{equation*}
Here the operator $\mathcal{L}_\psi:C^\infty(X)\rightarrow C^\infty(X)$ is given by the derivative of a function $f$ in the direction of the vector field $X$, i.e.
\begin{equation*}
\mathcal{L}_\Psi f:=\left. \frac{d}{dt}\right|_{t=0}\widehat{\psi}_t(f),
\end{equation*}
where $\psi^t$ denotes the flow of $\Psi$. Moreover, one can show that the map is a Lie-algebra anti-homomorphism:
\begin{equation*}
\mathcal{L}_{[\Psi,\Phi]}=\mathcal{L}_\Phi\mathcal{L}_\Psi-\mathcal{L}_\Psi\mathcal{L}_\Phi
\end{equation*}
for any vector fields $\Psi,\Phi\in \mathrm{Vect}(X)$.\\

\noindent
Hence, using the structure preserving condition \eqref{identity} and differentiating Diagram \eqref{diagram3} at the identity we obtain the new diagram:
\begin{equation*}
\begin{tikzcd}
C^\infty(X) \arrow[r,"\mathcal{L}_\Psi"] \arrow[d,"\pi^n"'] & C^\infty(X) \arrow[d,"\pi^n"] \\
A_N\arrow[r,"\mathcal{L}_{\Psi_N}"] & A_N
\end{tikzcd}
\end{equation*}
where $\mathcal{L}_\psi$ and $\mathcal{L}_{\psi_N}$ are derivations.
\begin{prop}
A structure preserving discretization of the diffeomorphism group $\mathrm{Diff}(X)$ induces a structure preserving discretization of the Lie algebra $(\mathrm{Vect}(X),\mathcal{L}_\psi)$.
\begin{itemize}
\item[i)]$C^\infty(X)$ is merely a vector space, $\mathfrak{D}(\mathrm{Vect}(X))=\mathfrak{gl}_N(\mathbb{C}).$
\item[ii)]$(C^\infty(X),\left\langle \cdot,\cdot \right\rangle )$ is a prehilbert space, then $\mathfrak{D}(\mathrm{Vect}_0^{+}(X))=\mathfrak{so}_N(\mathbb{C}).$
\item[iii)]$(C^\infty(X),1)$ is a vector space with a unit element (i.e. containing its scalar field), then $\mathfrak{D}(\mathrm{Vect}_{\mathrm{m}}(X))=\mathfrak{gl}_N(\mathbb{C}).$
\end{itemize}
\end{prop}
\subsection{Berezin-Toeplitz Quantization}
The Berezin-Toeplitz quantization bridges the gap between classical and quantum mechanics by transitioning from the Poisson algebra of smooth functions on a Kähler manifold to a finite-dimensional, noncommutative matrix algebra. This quantization method maintains essential geometric and algebraic properties, making it an ideal framework for discretizing continuous structures while ensuring structural integrity. In this section, the detailed process of Berezin-Toeplitz quantization is presented, and we demonstrate its role in the discretization of compact manifolds. We prove from first principle that it is a structure preserving method. It also provides a significant example to understand how continuous spaces are represented discretely in noncommutative settings.\\

\noindent
The natural setting of the Berezin-Toeplitz quantization is a compact Kähler manifold $(M,\omega)$. Whereas the basic objects can be defined on an arbitrary Kähler manifold, some proofs rely on compactness to avoid technicalities. For this reason, we restrict ourselves to the compact setting.\\

\noindent
A Kähler manifold $(M,\omega)$ is called \textit{quantizable} if there exists an associated quantum line bundle $(L,h,\nabla)$, where $L$ is a holomorphic line bundle $L$ over $M$ (i.e. locally a one-dimensional vector space), $h$ a Hermitian metric on $L$, and $\nabla$ a connection compatible with the metric $h$ and the complex structure.\\

\noindent
We consider the space of holomorphic sections $\mathrm{\Gamma}_{hol}(M,L)$ with value in the line bundle $L$. The bundle $L$ will be the stage for a finite dimensional quantization. For this, one needs to define the appropriate Hilbert space from $\mathrm{\Gamma}_{hol}(M,L)$. the scalar product on the space of sections and the norm are given by
\begin{equation*}
\left\langle \varphi,\psi \right\rangle :=\int_Mh(\varphi,\psi)\mathrm{\Omega}, \quad \|\varphi\|:=\sqrt{\left\langle \varphi,\varphi \right\rangle },
\end{equation*}
where $\mathrm{\Omega}=\frac{1}{n!}\omega^{\wedge n}$ is the Liouville form, used as a volume form on $M$. Let $\mathrm{L}^2(M,L)$ be the $\mathrm{L}^2$-completion of $\mathrm{\Gamma}_\infty(M,L)$, and $\mathrm{\Gamma}_{hol}(M,L)$ be its finite dimensional (due to compactness) subspace of holomorphic sections. Let $\mathrm{\Pi}^{(1)}:\mathrm{L}^2(M,L)\rightarrow \mathrm{\Gamma}_{hol}(M,L)$ be the projection. For $f\in C^\infty(M)$, the \textit{Toeplitz operator} $T^{(1)}_f$ is defined to be 
\begin{equation}
T^{(1)}_f:=\mathrm{\Pi}^{(1)}(f\cdot): \quad \mathrm{\Gamma}_{hol}(M,L) \rightarrow \mathrm{\Gamma}_{hol}(M,L)
\end{equation}
The linear map 
\begin{equation}
T:C^\infty(M)\rightarrow \mathrm{End}\left(\mathrm{\Gamma}_{hol}(M,L) \right), \qquad f,\mapsto T_f
\end{equation}
is the \textit{Berezin-Toeplitz quantization map}. In general, because
\begin{equation}
T_fT_g\neq T_{fg}
\end{equation}
the quantization map is neither a Lie algebra homomorphism nor an associative algebra homomorphism. \\

\noindent
Moreover, the Berezin-Toeplitz quantization map from the commutative algebra of functions to a noncommutative finite-dimensional matrix algebra. The finite dimensionality implies the loss of a lot of classical information. In order to recover this information, one should consider not just the bundle $(L,\nabla, h)$ alone but all its tensor powers $(L^{\otimes m},\nabla^{(m)},h^{(m)})$ and apply the construction for every $m\in \mathbb{N}_0$. In this way, the direct sum of the spaces of holomorphic sections of $(L^{\otimes m},\nabla^{(m)},h^{(m)})$ gets identified with a Hilbert subspace $\mathcal{H}$ called generalized Hardy space.\\

\noindent
Thus, we have the following decomposition:
\begin{equation*}
\mathcal{H}=\bigoplus_m(\mathcal{H}_m,\left\langle \cdot,\cdot \right\rangle_m),\qquad \mathcal{H}_m=\mathrm{\Gamma}_{hol}(M,L^m).
\end{equation*}
where we have defined the modified measure 
\begin{equation}
\left\langle \varphi,\psi \right\rangle_m :=\int_Mh(\varphi,\psi)\mathrm{\Omega}^{(m)}_\epsilon,\qquad \text{where}\quad\mathrm{\Omega}^{(m)}_\epsilon:=\epsilon^{(m)}(x)\mathrm{\Omega}(x)
\end{equation}
with $\epsilon^{(m)}$ the \textit{Rawnsley's epsilon function}.
\begin{defn}[Toeplitz operator]
For $f\in C^\infty(M)$ the Toeplitz operator $T_f^{(m)}$ (of level $m$) is defined by
\begin{equation*}
T^{(m)}_f:=\mathrm{\Pi}^{(m)}f\mathrm{\Pi}^{(m)}: \mathrm{\Gamma}_{hol}(M,L^m)\rightarrow \mathrm{\Gamma}_{hol}(M,L^m)
\end{equation*}
\end{defn}
\noindent
This infinite family should in some sense \textit{approximate} the algebra $C^\infty(M)$; one obtains a family of matrix algebras and maps
\begin{equation*}
T^{(m)}:C^\infty(M) \rightarrow \mathrm{End}(\mathrm{\Gamma}_{hol}(M,L^m)), \qquad f\mapsto T_f^{(m)}=\mathrm{\Pi}^{(m)}(f\cdot), m\in \mathbb{N}_0
\end{equation*}
\begin{prop}[\cite{schlichenmaier_berezin-toeplitz_2010}]
\label{surj}
The Toeplitz map 
\begin{equation*}
C^\infty(M)\rightarrow End(\mathcal{H}_m), \quad f\mapsto T^{(m)}(f)
\end{equation*}
is surjective.
\end{prop}
\begin{defn}[Berezin-Toeplitz quantization map]
The Berezin-Toeplitz quantization map is the map
\begin{equation*}
\mathrm{\Lambda} : C^\infty(M)\rightarrow \prod_{m\in \mathbb{N}_0}\text{End}\left( \mathrm{\Gamma}_{hol}\left( M,L^{(m)}\right) \right) , \quad f\rightarrow T_f:=\left( T_f^{(m)}\right)_{m\in \mathbb{N}_0}
\end{equation*}
\end{defn}
\noindent
In other words, for $f\in C^\infty(M)$ we have the decomposition
\begin{equation}
T_f=\prod_{m=0}^\infty T_f^{(m)}
\end{equation}
where $T_f^{(m)}$ are exactly the restriction of $T_f$ to $\mathcal{H}^{(m)}$.
\begin{thm}[Bordemann, Meinrenken, Schlichenmaier, 
\cite{schlichenmaier_berezin-toeplitz_2010}]~\\
\label{BMS}
(a) For every $f\in C^\infty(M)$ there exists a $C>0$ such that
\begin{equation}
\|f\|_\infty -\frac{C}{m}\leq \|T^{(m)}_f\|\leq \|f\|_\infty
\end{equation}
In particular, $\lim_{m\rightarrow +\infty}\|T_f^{(m)}\|=\|f\|_\infty$.\\
(b) For every, $f,g\in C^\infty(M)$
\begin{equation}
\|im[T_f^{(m)},T_g^{(m)}]-T^{(m)}_{\left\lbrace f,g\right\rbrace }\|=O(m^{-1})
\end{equation}
\end{thm}
\begin{defn}
The \textit{covariant Berezin symbol} $\sigma^{(m)}$ (of level $m$) is defined by 
\begin{equation*}
\sigma^{(m)}:\mathrm{End}(\mathrm{\Gamma}_{hol}(M,L^{(m)}))\rightarrow C^\infty(M), \quad A\mapsto x\mapsto \sigma^{(m)}(A)(x):=\frac{\left\langle e_q^{(m)},Ae_q^{(m)} \right\rangle }{\left\langle e_q^{(m)},e_q^{(m)} \right\rangle }
\end{equation*}
\end{defn}
\noindent
We introduce on $\mathrm{End}(\mathrm{\Gamma}_{hol}(M,L^{(m)}))$ the Hilbert-Schmidt norm 
\begin{equation*}
<A,C>_{HS} =\mathrm{Tr}(A^*C)
\end{equation*}
\begin{prop}[\cite{schlichenmaier_berezin-toeplitz_2010}] The following statements hold:\\
i) The Toeplitz map $f\mapsto T_f^{(m)}$ and the covariant symbol map $A\mapsto \sigma^{(m)}(A)$ are adjoint
\begin{equation*}
\left\langle A,T_f^{(m)}\right\rangle_{HS}=\left\langle \sigma^{(m)}(A),f\right\rangle_m 
\end{equation*}
ii) The covariant symbol map $\sigma^{(m)}$ is injective.\\
iii) The \textit{Berezin transform} map defined
\begin{equation*}
C^\infty(M)\rightarrow C^\infty(M), \qquad f\mapsto I^{(m)}(f):=\sigma^{(m)}\circ T_f^{(m)}
\end{equation*}
satisfies the limit 
\begin{equation*}
\|f-I^{(m)}(f)\|_\infty \rightarrow 0 \quad \text{as $m\rightarrow \infty$}.
\end{equation*}
\end{prop}
\noindent
In this section, we are going to draw a bridge between the classical problem of \textit{structure preserving discretization} and the Berezin-Toeplitz quantization. We are going to show that the Berezin-Toeplitz quantization is in fact a structure preserving discretization of the Poisson algebra $(C^\infty(M), \left\lbrace \cdot, \cdot \right\rbrace )$, when $M$ is a compact Kähler manifold. \\

\noindent
Now that we have presented the background material on the Berezin-Toeplitz quantization, we can be described it by the following diagram.
\begin{equation*}
\begin{tikzcd}
C^\infty(M)\otimes C^\infty(M) \arrow[r,"{\left\lbrace \cdot, \cdot \right\rbrace} "] \arrow[d,"T^{(m)}"] & C^\infty(M) \arrow[d,"T^{(m)}"] \\
M_m(\mathbb{C})\otimes M_m(\mathbb{C}) \arrow[r,"{\left[ \cdot, \cdot\right] } "] & M_m(\mathbb{C}) 
\end{tikzcd}
\end{equation*}
Indeed, we prove that the Berezin-Toeplitz quantization satisfies the axiom of a structure preserving discretization of the pair $(C^\infty(M), \left\lbrace \cdot, \cdot \right\rbrace )$.
\begin{prop}
The Berezin-Toeplitz quantization of the Poisson algebra of smooth functions $(C^\infty(M), \left\lbrace \cdot, \cdot \right\rbrace )$ is a structure preserving discretization.
\end{prop}
\begin{proof}
We start by proving that the quantization is a discretization method in the sense of Definition \eqref{def1}. The Toeplitz maps $T^{(m)}$ are contractive by Theorem \eqref{BMS} and surjective by Proposition \eqref{surj}. Moreover, we have the limit:
\begin{equation*}
\lim_{m\rightarrow \infty}\|T^{(m)}(f)\|=\|f\|
\end{equation*}
again following from Theorem \eqref{BMS}. In addition, to any projection $T^{(m)}$, we can associate an injective section map $\sigma^{(m)}$. Therefore, the Berezin-Toeplitz quantization is a discretization.\\
Now, if we consider the pair $(C^\infty(M),\left\lbrace \cdot, \cdot\right\rbrace )$ as an object in the category of Lie algebra, then the pair $(M_m(\mathbb{C}),[\cdot, \cdot])$ is a structure preserving discretization since :
\begin{itemize}
\item[i)] $(M_m(\mathbb{C}),[\cdot, \cdot])$ is a Lie algebra
\item[ii)]  $\|T^{(m)}\left\lbrace f,g \right\rbrace - [T^{(m)}(f),T^{(m)}(g)]\|\rightarrow 0$ as $m\rightarrow \infty$.
\end{itemize}
Therefore, the Berezin-Toeplitz quantization provides a structure preserving discretization of the Lie algebra structure on $C^\infty(M)$.
\end{proof}
\section{Noncommutative Laplacian as a discrete Laplacian}
\label{Sect4}
\noindent
The Laplacian, which traditionally arises in the context of smooth Riemannian manifolds, is critical for understanding various physical phenomena, such as heat diffusion and wave propagation. However, when transitioning to discrete spaces, particularly in curved spaces, the formulation of the Laplacian requires a novel approach.\\

\noindent
This section focuses on how the Berezin-Toeplitz quantization framework can be used to discretize the classical Laplacian, thereby creating a discrete version that operates on matrix algebras rather than smooth functions. The key idea is to show that the noncommutative Laplacian emerges naturally, following the axioms, as a discrete counterpart of the classical Laplace operator on compact Kähler manifolds. By leveraging the structure-preserving properties of the Berezin-Toeplitz quantization, the resulting discrete Laplacian retains important geometric and analytic features, offering a powerful tool for studying differential equations in noncommutative settings.\\

\noindent
Consider again a compact Kähler manifold $(M,\omega)$. The Nash embedding theorem states that any Riemannian manifold can be isometrically embedded in the Euclidean space $\mathbb{R}^d$ for sufficiently large $d$. Thus, for a closed Poisson manifold $(M,\pi )$ with a metric $g$, there exists an isometry embedding 
\begin{equation}
X:M\rightarrow \mathbb{R}^d
\end{equation}
for sufficiently large $d$. We denote the embedding coordinate functions as $\left\lbrace X^k \right\rbrace_{k=1,\dots, d}$. The word \textit{isometric} means that the induced metric of the embedding is equal to the intrinsic metric $g$ on $M$: 
\begin{equation}
g_{ij}(x)=g_x(\partial_i,\partial_j)=\left\langle d_xX(\partial_i),d_xX(\partial_j) \right\rangle = \sum_{k=1}^d(\partial_iX^k)(\partial_jX^k) 
\label{coord}
\end{equation}
Using the smooth bivector $\pi\in \mathfrak{X}(M)$, we can associate to any function $f$ a derivation through the map defined by:
\begin{equation}
\xi_f:C^\infty(M)\rightarrow \mathrm{Der}(M), \quad \xi_f(g)=\pi(df\wedge dg)
\end{equation}
We can now define the following family of self-adjoint differential operators
\begin{equation}
\partial_k:C^\infty(M)\rightarrow C^\infty(M), \quad g\mapsto \partial_k(g):=\xi_{X^k}(g) =\left\lbrace X^k,g\right\rbrace 
\end{equation}
for every $k=1,\dots ,d$.
\begin{rmk}
One can then construct $\mathfrak{g}:=\mathrm{Lie}(<\partial_k>)$ the Lie algebra generated by the set of vector fields $(\partial_k)$ and define the so-called quantized Weil-algebra by
\begin{equation}
\mathcal{W}(\mathfrak{g})=\mathrm{U}(\mathfrak{g})\otimes \mathrm{Cl}(\mathfrak{g})
\end{equation}
where $\mathrm{U}(\mathfrak{g})$, respectively $\mathrm{Cl}(\mathfrak{g})$, is the enveloping Lie algebra, respectively the Clifford algebra, associated to $\mathfrak{g}$. One can then define a first order elliptic differential operator called the Dirac operator
\begin{equation}
D=\sum_k\partial_ke^k + \sum_{ijk}\gamma_{ijk}e^ie^je^k
\end{equation}
which plays a central role in the Attiyah-Singer theory. In this work however, we will simply consider the first order operator $\mathcal{D}$ given by 
\begin{equation}
\mathcal{D}= \sum_k \left\lbrace X^k,\cdot \right\rbrace .
\end{equation}
\end{rmk}
\noindent
We also define the Bochner-Laplace operator for the metric $g$ on $M$ defined by 
\begin{equation}
\mathrm{\Delta}\varphi := -g^{ij}\nabla_i\nabla_j\varphi. 
\end{equation}
using the covariant derivatives. There is a relation between the operator $D$ and the Laplacian given by
\begin{equation}
D^2=\mathrm{\Delta} +\text{lower order terms}
\end{equation}

\vspace{8pt}
\noindent
Therefore, we will derive structure preserving discretization of the differential algebra $(C^\infty(M),\mathcal{D})$ which will induce a discretization of the Laplacian.
\begin{thm}
The Berezin-Toeplitz quantization induces a structure preserving discretization of the differential algebra $(C^\infty(M),\mathcal{D})$ given by the following commuting diagram
\begin{equation}
\begin{tikzcd}
C^\infty(M) \arrow[r,"\mathcal{D}"] \arrow[d,"T^{(m)}"'] & C^\infty(M) \arrow[d,"T^{(m)}"] \\
\mathrm{End}(\mathcal{H}_m)\arrow[r,"d_m"] & \mathrm{End}(\mathcal{H}_m)
\end{tikzcd}
\end{equation}
where the differential $d_m$ is given by the following commutator:
\begin{equation}
A\mapsto d_m(A)=[D_m,A], \qquad D_m:=\sum_{k=1}^d\partial_m^k:=\sum_{k=1}^d \mathrm{\Pi}^{(m)}X_k\mathrm{\Pi}^{(m)}
\end{equation}
\end{thm}
\begin{proof}
We consider the pair $(C^\infty(M),\mathcal{D} )$ as an object in the category of differentiable algebra, then the pair $(M_m(\mathbb{C}),d_m)$ is a structure preserving discretization since :
\begin{itemize}
\item[i)] $(M_m(\mathbb{C}),[\cdot, \cdot])$ is a differential algebra.
\item[ii)]  $\|T^{(m)}\left\lbrace X^k,f \right\rbrace - [T^{(m)}(X^k),T^{(m)}(f)]\|\rightarrow 0$ as $m\rightarrow \infty$.
\end{itemize}
where the last statement follows from the quantization of the Poisson bracket.
\end{proof}
\begin{cor}
\label{cor2}
The Berezin-Toeplitz discretization $\mathfrak{D}(C^\infty(M),\mathcal{D})$ defines a discretization  
\begin{equation}
\mathrm{\Delta}_m: \mathrm{End}(\mathcal{H}_m)\rightarrow \mathrm{End}(\mathcal{H}_m),\qquad A\mapsto \mathrm{\Delta}_m(A)=\sum_{k=1}^d\left[ \partial_m^k\left[ \partial_m^k,A\right] \right] 
\end{equation}
of the Bochner-Laplace operator $\mathrm{\Delta}$ such that 
\begin{equation}
D_m^2=\mathrm{\Delta}_m+\text{lower order terms}
\end{equation}
\end{cor}
\begin{proof}
Recall the two identities that determine the quantization of the Poisson algebra:
\begin{align*}
 T_m(fg)=T_m(f)T_m(g) + O\left(m^{-1}\right) \quad \text{and} \quad  m\left[ T_m(f),T_m(g)\right] =T_m(\left\lbrace f,g \right\rbrace )+O(m^{-1}).
\end{align*}
Then, applying successively the identities to the iterated Poisson brackets defining the Laplace operator, we get:
\begin{align*}
T_m\left( \left\lbrace X^k , \left\lbrace X^k,\varphi \right\rbrace \right\rbrace\right) &= m\left[ T_m(X^k),T_m\left( \left\lbrace X^k,\varphi \right\rbrace \right)\right] + O(m^{-1}) \\
&= m^2\left[ T_m(X^k),\left[ T_m(X^k),T_m(\varphi)\right] + O(m^{-2})1\right] + O(m^{-1})\\
&=m^2\left[ T_m(X^k),\left[ T_m(X^k),T_m(\varphi)\right] \right] + \left[ T_m(X^k), O(1)1\right] + O(m^{-1})\\
&=\Delta_m(\varphi) + O(m^{-1})
\end{align*}
The statement then follows by take the norm on both sides of the last equality.
\end{proof}
\begin{cor}
The \textit{discrete Bochner-Laplace} operator
\begin{equation}
\widetilde{\mathrm{\Delta}}_m:C^\infty(M)\rightarrow C^\infty(M), \qquad \widetilde{\mathrm{\Delta}}_m(f)=\sigma^{(m)}\circ \mathrm{\Delta}_m\circ T^{(m)}
\end{equation}
is a self-adjoint operator on $(\mathcal{H}_m,<\cdot,\cdot>_m)$ with discrete spectrum.
\end{cor}
\begin{proof}(Recall that for $f\in C^\infty$, the operator $M_f$ is multiplication with $\tau^*f$).\\
Using the spectral theorem applied to $\mathrm{\Delta}_m$ we get the decomposition of $\widetilde{\mathrm{\Delta}}_m$ as
\begin{equation}
\widetilde{\mathrm{\Delta}}_m:\mathcal{H}_m\rightarrow \mathcal{H}_m,\qquad \widetilde{\mathrm{\Delta}}_m(f)=\sum_{i=1}^m\lambda_i\sigma^{(m)}P_{\lambda_i}T^{(m)}
\end{equation}
where $(\lambda_i,P_{\lambda_i})$ are the pairs of eigenvalues and associated spectral projections.\\
The operator $\sigma^{(m)}P_{\lambda}T^{(m)}$ is a self-adjoint projection for every $\lambda\in \sigma(\mathrm{\Delta}_m)$.
\end{proof}
\noindent
An important question that remains to address is the convergence of the discrete Bochner-Laplacian to the continuous counterpart. In order to do so, we need to show, following Definition \eqref{conv}, the limit 
\begin{equation*}
\lim_{m\rightarrow \infty}\|\mathrm{\Delta}\varphi-\sigma^{(m)}\circ \mathrm{\Delta}_m\circ T^{(m)}\|=0
\end{equation*}
We will need the following definition to tackle the convergence of unbounded linear operators.
\begin{defn}[Strong graph limit]
Let $A_n$ be a sequence of operators on a Hilbert space $\mathcal{H}$. We say that $(\psi,\varphi)\in \mathcal{H}\times \mathcal{H}$ is in the \textit{strong graph limit} of $A_n$ if we can find $\psi_n\in D(A_n)$ so that $\psi_n\rightarrow \psi$, $A_n\psi_n\rightarrow \varphi$. We denote the set of pairs in the strong graph limit by $\mathrm{\Gamma}^s$. If $\mathrm{\Gamma}^s$ is the graph of an operator $A$ we say that $A$ is the \textit{strong graph limit} of $A_n$ and write $A=\mathrm{sgr-lim}A_n$.
\end{defn}
\begin{thm}
\label{thm-lp}
The discrete Bochner-Laplace operator $\widetilde{\mathrm{\Delta}}_m$ converges in the strong-graph limit sense to the Bochner-Lapalce operator $\mathrm{\Delta}$:
\begin{equation*}
\mathrm{sgr}-\lim \widetilde{\mathrm{\Delta}}_m = \mathrm{\Delta}
\end{equation*}
\end{thm}
\noindent
In order to prove the theorem, we need the following lemma relating the Laplacian to the Poisson bracket.
\begin{lem}
Let $(M,g,\omega)$ be a Kähler manifold, with its induced Poisson structure $(M,\pi)$, such that $(M,g)$ which admits an isometric embedding $X:M\rightarrow \mathbb{R}^d$. Then, the Laplacian operator can be expressed as:
\begin{equation}
\mathrm{\Delta}_g\varphi = -\sum_{k=1}^d\left\lbrace X^k , \left\lbrace X^k,\varphi \right\rbrace \right\rbrace
\end{equation}
\end{lem}
\begin{proof}
We start by expressing the Poisson bracket using the bivector $\pi$ and the connection operator:
\begin{align*}
-\left\lbrace X^k , \left\lbrace X^k,\varphi \right\rbrace \right\rbrace &=- \pi\left(dX^k\wedge \pi\left( dX^k\wedge d \varphi\right) \right) =-\pi^{ij}\pi^{rs}\partial_iX^k\nabla_j\left(\partial_rX^k \nabla_s\varphi \right)
\end{align*}
Using the product rule, we obtain:
\begin{align*}
&-\pi^{ij}\pi^{rs}\partial_iX^k\nabla_j\left(\partial_rX^k \nabla_s\varphi \right)=-\pi^{ij}\pi^{rs}\partial_iX^k\left[ \nabla_j\left(\partial_rX^k\right)\nabla_s\varphi  +  \partial_rX^k\nabla_j\nabla_s\varphi \right]\\
&=-\pi^{ij}\pi^{rs}\left[\nabla_j\left(\partial_iX^k\partial_rX^k\right)\nabla_s\varphi  -  (\nabla_j\partial_iX^k)(\partial_rX^k)\nabla_s\varphi + \partial_iX^k\partial_rX^k\nabla_j\nabla_s\varphi\right] .
\end{align*}
We have used the fact that $X$ are coordinate embedding and satisfy Equation \eqref{coord}. We use now that $\nabla_j\left(g_{ir}\right)=0$ (metric compatible connection) and
\begin{equation*}
\pi^{ij}\nabla_j\partial_iX^k=\pi^{ij}\left(\partial_j\partial_i X^k - \Gamma_{ij}^\ell\partial_\ell X^k \right)= \pi^{ij}\partial_j\partial_i X^k - \pi^{ij}\Gamma_{ij}^\ell\partial_\ell X^k=0
\end{equation*}
 since $\pi^{ij}=-\pi^{ji}$ and also $\Gamma^k_{ij}=\Gamma^k_{ji}$. We deduce then that
\begin{align*}
-\left\lbrace X^k , \left\lbrace X^k,\varphi \right\rbrace \right\rbrace=-\pi^{ij}\pi^{rs}\left[\nabla_j\left(g_{ir}\right)\nabla_s\varphi  -  (\nabla_j\partial_iX^k)(\partial_rX^k)\nabla_s\varphi + g_{ir}\nabla_j\nabla_s\varphi\right]
\end{align*}
 Important, here we have used $\nabla \pi =0$ which follows from the general properties of the Kähler structure $\nabla g =\nabla J=\nabla \omega =0$. Therefore, we are left with
\begin{equation*}
-\left\lbrace X^k , \left\lbrace X^k,\varphi \right\rbrace \right\rbrace = -\pi^{ij}\pi^{rs} g_{ir}\left(\nabla_j\nabla_s\varphi\right)=-g^{js} \left(\nabla_j\nabla_s\varphi\right).
\end{equation*}
\end{proof}
\begin{lem}
\label{lem-lap2}
Let $(M,g,\omega)$ be a Kähler manifold, with its induced Poisson structure $(M,\pi)$, such that $(M,g)$ which admits an isometric embedding $X:M\rightarrow \mathbb{R}^d$. Then, the bi-Laplacian operator can be expressed as:
\begin{equation*}
\sum_{k}\left\lbrace X^k, \mathrm{\Delta}\left\lbrace X^k ,\varphi \right\rbrace \right\rbrace =\mathrm{\Delta}^2(\varphi)
\end{equation*}
for every $\varphi\in C^\infty(M)$.
\end{lem}
\begin{proof}
We apply the Laplacian to the Poisson bracket with $X^k$
\begin{equation*}
\mathrm{\Delta}\left\lbrace X^k ,\varphi \right\rbrace = \left\lbrace \mathrm{\Delta}X^k ,\varphi \right\rbrace + \left\lbrace \nabla_i X^k ,\nabla_i\varphi \right\rbrace  +\left\lbrace X^k ,\mathrm{\Delta}\varphi \right\rbrace = \left\lbrace X^k ,\mathrm{\Delta}\varphi \right\rbrace
\end{equation*}
using the fact that $\pi^{ij}\nabla_j\partial_iX^k=0$.
\end{proof}
\begin{lem}
\label{Lem-cont}
The section maps $\sigma^{(m)}$ are contractives i.e.
\begin{equation*}
\|\sigma^{(m)}(A)\|\leq \|A\|
\end{equation*}
for every $A\in M_m(\mathbb{C})$.
\end{lem}
\begin{lem}
(a) For every $f\in C^\infty(M)$, we have the first estimate
\begin{equation}
I^{(m)}(f)=f+\frac{1}{m}\mathrm{\Delta}(f)+O(m^{-2})
\label{eq-I}
\end{equation}
(b) For every, $f,g\in C^\infty(M)$, we have the second estimate
\begin{equation}
\mathrm{i}m[T^{(m)}(\varphi),T^{(m)}(\psi)]=T^{(m)}\left( \left\lbrace f,g \right\rbrace \right)+O(m^{-1})T^{(m)}(\mathrm{\Delta}\left\lbrace f,g \right\rbrace )
\end{equation}
\end{lem}
\begin{proof} Theorem \ref{thm-lp}\\
We start by applying a triangular inequality
\begin{align*}
&\|\mathrm{\Delta}(\varphi)-\sigma^{(m)}\circ\mathrm{\Delta}_m\circ T^{(m)}(\varphi)\|_{\infty} \\
&\leq \|\mathrm{\Delta}(\varphi)-\sigma^{(m)}\circ T^{(m)}(\mathrm{\Delta}(\varphi))\|_{\infty}
+ \|\sigma^{(m)}\circ T^{(m)}(\mathrm{\Delta}(\varphi)) - \sigma^{(m)}\circ\mathrm{\Delta}_m\circ T^{(m)}(\varphi)\|_{\infty} 
\end{align*}
and then we apply Lemma \ref{Lem-cont} to bound the second term in the right-hand-side:
\begin{align*}
\|\mathrm{\Delta}(\varphi)-&\sigma^{(m)}\circ\mathrm{\Delta}_m\circ T^{(m)}(\varphi)\|_{\infty} \\
&\leq \|\mathrm{\Delta}(\varphi)-\sigma^{(m)}\circ T^{(m)}(\mathrm{\Delta}(\varphi))\|_{\infty} + \| T^{(m)}(\mathrm{\Delta}(\varphi)) - \mathrm{\Delta}_m\circ T^{(m)}(\varphi)\|_{\infty}
\end{align*}
We then combine \eqref{eq-I} and Lemma \ref{lem-lap2} to bound the first term in the right-hand-side and Corollary \ref{cor2} to bound the second term to get
\begin{equation*}
\|\mathrm{\Delta}(\varphi)-\sigma^{(m)}\circ\mathrm{\Delta}_m\circ T^{(m)}(\varphi)\|_{\infty} \leq a_m\|\mathrm{\Delta}^2(\varphi)\|_{\infty}
\end{equation*}
Since we work on a compact smooth manifold, we have one hand that 
\begin{equation*}
\|\mathrm{\Delta}(\varphi)-\sigma^{(m)}\circ\mathrm{\Delta}_m\circ T^{(m)}(\varphi)\|_{2}\leq c_0\|\mathrm{\Delta}(\varphi)-\sigma^{(m)}\circ\mathrm{\Delta}_m\circ T^{(m)}(\varphi)\|_{\infty}
\end{equation*}
and on the other using Sobolev embedding we have
\begin{equation*}
\|\mathrm{\Delta}^2(\varphi)\|_{\infty}\leq c_1\|\varphi\|_{H^k}
\end{equation*}
In conclusion, we get 
\begin{equation*}
\|\mathrm{\Delta}(\varphi)-\sigma^{(m)}\circ\mathrm{\Delta}_m\circ T^{(m)}(\varphi)\|_{2} \leq a_m\|\varphi\|_{H^k}
\end{equation*}
\end{proof}
\begin{rmk}
It is important to notice that the structure preserving property obtained in Corollary \ref{cor2} plays a central role in the proof of convergence; hence the last theorem is an example of results showing how structure preservation might imply convergence.
\end{rmk}
\noindent
We can also answer the question about the convergence of the spectrum of the Laplacian, which is crucial in many applications.
\begin{thm}
Let $A_n$ $(n\in \mathbb{N})$ and $A$ be (unbounded) self-adjoint operators and assume that $D(A)=D(A_n)$. Assume, furthermore, that there are null-sequences $(a_n)$ and $(b_n)$ from $\mathbb{R}$ for which
\begin{equation*}
\|(A-A_n)f\|\leq a_n\|f\|+b_n\|Af\| \quad \text{for all $f\in D(A)$}.
\end{equation*}
Then $\sigma(A)=\lim_{n\rightarrow \infty}\sigma(A_n)$.
\end{thm}
\begin{cor}
The spectrum discrete Bochner-Laplace operator $\widetilde{\mathrm{\Delta}}_m$ converges in the spectrum of the true Bochner-Lapalce operator $\mathrm{\Delta}$ i.e. $\sigma(\mathrm{\Delta})=\lim_{n\rightarrow \infty}\sigma(\mathrm{\Delta}_n)$.
\end{cor}
\section{Example of the sphere $\mathbb{S}^2$}
\label{Sect5}
\noindent
We show in the concrete example of the sphere $(\mathbb{S}^2,\lbrace\cdot, \cdot\rbrace)$, seen as a compact Poisson manifold, how a structure-preserving discretization of the differential structure of the space $\mathbb{S}^2$ can be derived from first principles, including the discretization of the Laplacian operator. This is achieved following the set of axioms introduced in this work; the projection maps are provided by the Berezin-Toeplitz quantization. As a result, we recover the discrete Laplacian used in Zeitlin's model for Euler's equations on the sphere \cite{modin_eulerian_2023,modin_liepoisson_2020}; hence proving that such model is structure-preserving in the sense of the present work. In a follow-up work \cite{tageddine_Laplace_2025}, we extend this construction to compact surfaces with higher genus. \\

\noindent
Let us recall that there exists an isometric embedding:
\begin{equation}
\rho : \mathbb{S}^2\rightarrow \mathbb{R}^3, \qquad \rho(s):=(x,y,z)
\end{equation}
of the sphere onto $\mathbb{R}^3$, such that $x^2+y^2+z^2=1$. We consider the differential operator $\mathcal{D}$ defined by 
\begin{equation}
\mathcal{D}(f)=\lbrace x,f\rbrace + \lbrace y,f\rbrace + \lbrace z,f\rbrace = - \left(\partial_xf + \partial_yf + \partial_zf \right) 
\end{equation}
Finding a structure-preserving discretization of the differential algebra $(C^\infty(\mathbb{S}^2),\mathcal{D})$ is equivalent by our axiomatization to the existence of the following diagram
\begin{equation}
\begin{tikzcd}
C^\infty(\mathbb{S}^2) \arrow[r,"\mathcal{D}"] \arrow[d,"T^{(m)}"'] & C^\infty(\mathbb{S}^2) \arrow[d,"T^{(m)}"] \\
\mathrm{End}(\mathcal{H}_m)\arrow[r,"d_m"] & \mathrm{End}(\mathcal{H}_m)
\end{tikzcd}
\end{equation}
such that $(\mathrm{End}(\mathcal{H}_m),d_m)$ is a differential algebra and $\|T^{(m)}\circ\mathcal{D} - d_m\circ T^{(m)}\|\rightarrow 0$.\\

\noindent
If we rewrite the spherical polynomial the spherical harmonics $\left\lbrace  Y_{\ell m}(\theta, \varphi)\right\rbrace  $ as the space $\mathcal{P}_{\mathrm{h}}[x_1,x_2,x_3]$ of harmonic homogeneous polynomial in $3$ variables $x_1,x_2,x_3$, then $T^{(m)}$ can be taken as the projection onto  the subspace $\mathcal{P}^m_{\mathrm{h}}[x_1,x_2,x_3]$ of degree $m$ homogeneous polynomials. In this case, $T^{(m)}$ corresponds to the Berezin-Toeplitz quantization of $\mathbb{S}^2$ and it satisfies the condition of a structure preserving discretization.\\

\noindent
Following the procedure of the previous section, the continuous operators $(x,y,z)$ are replaced by $m$-dimensional matrices:
\begin{equation}
T^{(m)}(x)=X_m, \quad T^{(m)}(y)=Y_m, \quad T^{(m)}(z)=Z_m,
\end{equation}
and the Lie algebra generated by $(X_m,Y_m,Z_m)$ is the Spin algebra $\mathfrak{so}(m)$. Then, the differential operator $\mathcal{D}$ is discretized using the Lie bracket of matrices and gives the derivation $d_m$:
\begin{equation}
d_m(A)=[X_m,A]+[Y_m,A]+[Z_m,A]
\end{equation}
Similarly, the Laplace operator on the sphere is replaced by 
\begin{equation}
\Delta_m(A)=[X_m,[X_m,A]]+[Y_m,[Y_m,A]]+[Z_m,[Z_m,A]],
\end{equation}
which is identified as the non-commutative Laplacian.
\nocite{*}
\addcontentsline{toc}{section}{References}
\bibliographystyle{amsplain}
\bibliography{References.bib}

\providecommand{\bysame}{\leavevmode\hbox to3em{\hrulefill}\thinspace}
\providecommand{\MR}{\relax\ifhmode\unskip\space\fi MR }
\providecommand{\MRhref}[2]{%
  \href{http://www.ams.org/mathscinet-getitem?mr=#1}{#2}
}
\providecommand{\href}[2]{#2}
\begin{thebibliography}{10}

\bibitem{arnold_finite_2010}
D.~N. Arnold, R.~S. Falk, and R.~Winther, \emph{Finite element exterior
  calculus: from {Hodge} theory to numerical stability}, Bulletin of the
  American Mathematical Society \textbf{47} (2010), no.~2, 281--354 (en),
  arXiv: 0906.4325.

\bibitem{arveson_discretized_nodate}
William Arveson, \emph{{Discretized} {CCR} {Algebras}} (en).

\bibitem{bordemann_toeplitz_1994}
Martin Bordemann, Eckhard Meinrenken, and Martin Schlichenmaier, \emph{Toeplitz
  quantization of {Kähler} manifolds and gl({N}), {$N\rightarrow \infty$}
  limits}, Communications in Mathematical Physics \textbf{165} (1994), no.~2,
  281--296 (en).

\bibitem{christiansen_topics_2011}
S.~H. Christiansen, H.~Z. Munthe-Kaas, and B.~Owren, \emph{Topics in
  structure-preserving discretization}, Acta Numerica \textbf{20} (2011),
  1--119 (en).

\bibitem{christiansen_construction_2008}
Snorre~H. Christiansen, \emph{{A} {C}onstruction of {S}paces of {C}ompatible
  {D}ifferential {F}orms on {C}ellular {C}omplexes}, Mathematical Models and
  Methods in Applied Sciences \textbf{18} (2008), no.~05, 739--757 (en).

\bibitem{connes_noncommutative_1994}
A.~Connes, \emph{Noncommutative geometry}, Academic Press, San Diego, 1994
  (eng).

\bibitem{connes_noncommutative_2007}
A.~Connes and M.~Marcolli, \emph{Noncommutative {Geometry}, {Quantum} {Fields}
  and {Motives}}, Colloquium {Publications}, vol.~55, American Mathematical
  Society, Providence, Rhode Island, December 2007 (en).

\bibitem{connes_spectral_2021}
Alain Connes and Walter~D. van Suijlekom, \emph{Spectral truncations in
  noncommutative geometry and operator systems}, Communications in Mathematical
  Physics \textbf{383} (2021), no.~3, 2021--2067 (en), arXiv:2004.14115
  [hep-th].

\bibitem{desbrun_discrete_2005}
M.~Desbrun, A.~N. Hirani, M.~Leok, and J.~E. Marsden, \emph{Discrete {Exterior}
  {Calculus}}, arXiv:math/0508341 (2005) (en), arXiv: math/0508341.

\bibitem{gawlik_geometric_2011}
E.S. Gawlik, P.~Mullen, D.~Pavlov, J.E. Marsden, and M.~Desbrun,
  \emph{Geometric, variational discretization of continuum theories}, Physica
  D: Nonlinear Phenomena \textbf{240} (2011), no.~21, 1724--1760 (en).

\bibitem{hagen_c-algebras_2001}
R.~Hagen, S.~Roch, and B.~Silbermann, \emph{C*-algebras and numerical
  analysis}, Monographs and textbooks in pure and applied mathematics, no. 236,
  Marcel Dekker, New York, 2001.

\bibitem{hirani_discrete_2003}
A.~N. Hirani, \emph{Discrete {Exterior} {Calculus}}, Ph.D. thesis, Caltech,
  Pasadena, California, 2003.

\bibitem{lundervold_hopf_2010}
Alexander Lundervold and Hans Munthe-Kaas, \emph{Hopf algebras of formal
  diffeomorphisms and numerical integration on manifolds}, September 2010,
  arXiv:0905.0087 [math].

\bibitem{mclachlan_butcher_2017}
Robert~I. McLachlan, Klas Modin, Hans Munthe-Kaas, and Olivier Verdier,
  \emph{Butcher series: {A} story of rooted trees and numerical methods for
  evolution equations}, February 2017, arXiv:1512.00906 [math].

\bibitem{modin_eulerian_2023}
Klas Modin and Manolis Perrot, \emph{Eulerian and {Lagrangian} {Stability} in
  {Zeitlin}’s {Model} of {Hydrodynamics}}, Communications in Mathematical
  Physics \textbf{405} (2024), no.~8, 177 (en).

\bibitem{modin_liepoisson_2020}
Klas Modin and Milo Viviani, \emph{Lie–{Poisson} {Methods} for {Isospectral}
  {Flows}}, Foundations of Computational Mathematics \textbf{20} (2020), no.~4,
  889--921 (en).

\bibitem{sakai_c-algebras_1998}
S.~Sakai, \emph{C*-algebras and {W}*-algebras}, Classics in mathematics,
  Springer, Berlin ; New York, 1998.

\bibitem{sakai_operator_2008}
Shôichirô Sakai, \emph{Operator algebras in dynamical systems: the theory of
  unbounded derivations in {C}*-algebras}, Cambridge University Press,
  Cambridge, 2008 (eng), OCLC: 190967415.

\bibitem{schlichenmaier_berezin-toeplitz_2010}
Martin Schlichenmaier, \emph{Berezin-{Toeplitz} quantization for compact
  {Kaehler} manifolds. {A} {Review} of {Results}}, Advances in Mathematical
  Physics \textbf{2010} (2010), 1--38 (en), arXiv:1003.2523 [math-ph].

\bibitem{stetter_analysis_1973}
Hans~J. Stetter, \emph{Analysis of discretization methods for ordinary
  differential equations}, Springer, Berlin, 1973 (eng), OCLC: 681410907.

\bibitem{tageddine_noncommutative_2022}
Damien Tageddine and Jean-Christophe Nave, \emph{Noncommutative {Differential}
  {Geometry} on {Infinitesimal} {Spaces}}, accepted for publication in Annales Mathématiques du Québec (2025), arXiv:2209.12929.

\bibitem{tageddine_statistical_2023}
\bysame, \emph{Statistical {Fluctuation} of {Infinitesimal} {Spaces}}, passed first round of revisions in Annales Mathématiques du Québec (2025), arXiv:2304.10617.

\bibitem{tageddine_Laplace_2025}
\bysame, \emph{Noncommutative { Laplacian} and numerical approximation of Laplace-Beltrami spectrum of compact Riemann surfaces}, in preparation (2025).

\end{thebibliography}
\end{document}